\documentclass[11pt]{article}

\usepackage{amssymb,amscd,array}

\let\ssection=\section
\renewcommand{\section}{\setcounter{equation}{0}\ssection}

\def\d{\delta}

\def\om{\omega}
\def\r{\rho}
\def\a{\alpha}
\def\b{\beta}
\def\s{\sigma}
\def\vfi{\varphi}

\def\l{\lambda}
\def\m{\mu}
\def\n{\nabla}

\def\implies{\Rightarrow}

\newcommand{\bbR}{\mathbb{R}}

\newcommand{\Diff}{\mathrm{Diff}}

\newcommand{\cF}{{\mathcal{F}}}
\newcommand{\cD}{{\mathcal{D}}}

\newcommand{\Pol}{\mathrm{Pol}}

\newcommand{\Sl}{\mathrm{sl}}

\newcommand{\Vect}{\mathrm{Vect}}

\newcommand{\cqfd}{\hspace*{\fill}\rule{3mm}{3mm}}

\begin{document}



\def\d{\delta}
\def\g{\gamma}
\def\om{\omega}
\def\r{\rho}
\def\a{\alpha}
\def\b{\beta}
\def\s{\sigma}
\def\vfi{\varphi}
\def\l{\lambda}
\def\m{\mu}
\def\implies{\Rightarrow}

\oddsidemargin .1truein
\newtheorem{thm}{Theorem}[section]
\newtheorem{lem}[thm]{Lemma}
\newtheorem{cor}[thm]{Corollary}
\newtheorem{pro}[thm]{Proposition}
\newtheorem{ex}[thm]{Example}
\newtheorem{rmk}[thm]{Remark}
\newtheorem{defi}[thm]{Definition}

\title{Formula for the Projectively Invariant Quantization on Degree Three}
\author{Sofiane BOUARROUDJ\footnote{Research supported by the Japan Society for the Promotion of Science.}\\
{\footnotesize   
Department of Mathematics, Keio University, Faculty of Science \& Technology} \\
{\footnotesize 3-14-1, Hiyoshi, Kohoku-ku, Yokohama, 223-8522, Japan. }\\
{\footnotesize  Tel : 81 45 563 1141 - Fax : 81 45 563 5948}\\
{\footnotesize  mailto:sofbou@math.keio.ac.jp}
}
\date{}
\maketitle
\begin{abstract}
We give an explicit formula for the projectively invariant quantization map between the space of 
symbols of degree three and the space of third-order linear differential operators, both viewed as 
modules over the group of diffeomorphisms and the Lie algebra of vector fields on a manifold.
\end{abstract} 
\vskip 1cm
\begin{quotation}
\begin{center}
{\LARGE Une formule pour la quantification projectivement invariante en degr\'e trois}
\end{center}
\end{quotation}
\vskip 1cm

\begin{center}
{\bf R\'esum\'e}
\end{center}
\begin{quotation} {\small 
Nous donnerons une formule explicite pour la quantification projectivement invariante entre l'espace des 
symboles de degr\'es trois et l'espace des op\'erateurs diff\'erentiels lin\'eaires d'ordres trois, vus comme 
modules sur le groupe des diff\'eomorphismes et l'alg\`ebre de Lie des champs de vecteurs sur une vari\'et\'e 
diff\'erentiable.}
\end{quotation}
\section{Introduction}
Let $M$ be a manifold of dimension $n.$ Fix an affine connection $\nabla$ on $M.$ Denote by $\cF_{\l}(M)$ the 
space of $\l-$densities on $M$ (i.e. sections of the bundle $(\wedge^n T^* M)^{\otimes \l}$). This space admits 
naturally a structure of module over the group of diffeomorphisms $\Diff(M)$ and the Lie algebra of vector 
fields $\Vect(M).$ Consider $\cD_{\l,\mu}(M)$  the space of linear differential operators acting from 
$\cF_\l(M)$ to $\cF_\mu(M).$ This space is a module over $\Diff(M)$ and $\Vect(M)$ (see \cite{b,do,lo,l}). 
The action is given as follows: take $f\in\Diff(M)$ and $A\in \cD_{\l,\mu}(M)$ then
\begin{equation}
\label{op}
f^*A=f_{\mu}^*\circ A \circ {f_{\l}^*}^{-1},
\end{equation}
where $f_{\l}^*$ is the standard action of a diffeomorphism on $\cF_{\l}(M).$ 

Differentiating the action of the flow of a vector field, one gets the corresponding action of $\Vect(M).$

Denote by $\cD_{\l,\mu}^3(M)$ the space of third-order linear differential operators endowed with the structure 
of module (\ref{op}). The module $\cD_{\l,\mu}^3(M)$ is viewed as a submodule of $\cD_{\l,\mu}(M).$

Consider now $\Pol(T^*M)$ the space of functions on the cotangent bundle $T^*M,$ polynomials on the fibers. 
This space is naturally isomorphic to the space of symmetric contravariant tensor fields on $M.$ One can define a 
one-parameter family of $\Diff(M)-$modules by taking $\Pol_{\d}(T^*M):=\Pol(T^*M)\otimes \cF_\d (M).$ Let us give 
explicitly this action: take  $f\in\Diff(M)$ and $P\in \Pol_{\d}(T^*M)$ then 
\begin{equation}
\label{sy}
f^*_\d P=  f^* P\cdot  (J_f) ^{-\d},
\end{equation}
where $f^*$ is the natural action of a diffeomorphism on contravariant tensor fields, and $J_f$ is the Jacobian 
of $f.$ 

Differentiating the action of the flow of a vector field, one gets the corresponding action of $\Vect(M).$

Denote by $\Pol_{\d}^3(T^*M)$ the space of symbols of degree three endowed with the module structure given by 
(\ref{sy}). 

Suppose $M:=\bbR^n$ is endowed with a flat projective structure (coordinates change are projective transformations).  
In this case, Lecomte and Ovsienko in \cite{lo} construct a quantization map between the space $\Pol_{\d}(T^*\bbR^n)$ 
and the space $\cD_{\l,\mu}(\bbR^n)$, equivariant with respect to the action of the Lie algebra $\Sl_{n+1}(\bbR)\subset 
\Vect(\bbR^n).$ Consider now any manifold $M$ and fix an affine connection on it. It is interesting to ask if there 
exists a canonical quantization map associated to the given connection. On degree two, the author construct in \cite{b} 
a quantization map depending only on the projective class of the connection (see also \cite{do} for the conformal case). 
This approach generalizes Lecomte and Ovsienko's approach for the flat case. On higher order, the problem of existence of 
the projectively invariant quantization map is open.  
\section{Main result}
The main result of this note is 

\begin{thm}
\label{}
For $n>1,$ and for $\d\not= \frac{n+3}{n+1},\frac{n+4}{n+1},\frac{n+5}{n+1},$ there exits a quantization map 
$Q:\Pol_{\d}^3(T^*M)\rightarrow \cD_{\l,\mu}^3(M)$ 
given by 
\begin{equation}
\label{main}
\begin{array}{cl}
P^{ijk}\mapsto & P^{ijk}\n_i\n_j\n_k +\alpha \n_k P^{ijk}\,\n_i \n_j +\left (\beta_1 \n_i \n_j P^{ijk}+
\beta_2 P^{ijk} R_{ij}\right ) \n_k \\[3mm]
&+(\eta_1 \n_i\n_j\n_k P^{ijk}+ \eta_2 R_{ij} \n_k P^{ijk} +\eta_3 \n_i R_{jk}\,P^{ijk})
\end{array}
\end{equation}
where $R_{ij}$ are the components of the Ricci tensor of the connection $\nabla,$ the constants 
$\a,\b_1,\b_2,\eta_1, \eta_2,\eta_3,$ are given by 
$$
\begin{array}{ll}
\alpha= \displaystyle{\frac{6+3\l(1+n)}{4+(1-\d)(1+n)} },& \beta_1=\displaystyle{\frac{1+\l(n+1)} 
{3+(1-\d)(1+n)} \alpha }, \\[3mm]
 \beta_2=\displaystyle\frac{2+3\l(1+n)-(4+(1-\d)(1+n))\beta_1}{n-1},&
 \eta_1= \displaystyle{ \frac{\l (1+n)}{(6+3(1-\d)(1+n))}\beta_1},\\[3mm]
\eta_3=\displaystyle\frac{\l (1+n)-\eta_1 (4+(1-\d)(1+n))}{n-1},& \eta_2=\displaystyle\frac{\l (1+n)\alpha-
(10+3(1-\d)(1+n))\eta_1}{n-1},
\end{array}
$$
and have the following properties\\
(i) It depends only on the projective class of the connection $\nabla$ (see \cite{kn}).\\
(ii) If $M=\bbR^n$ is endowed with a flat projective structure the map (\ref{main}) is the unique map that 
preserves the principal symbols, equivariant with respect to the action of the Lie algebra  
$\Sl_{n+1}(\bbR)\subset \Vect(\bbR^n).$
\end{thm}
{\bf Proof.} Let us give an idea of the proof. Let $\tilde \nabla$ be another connection projectively 
equivalent to $\nabla.$ Denote by $Q^{\tilde \nabla}$ the quantization map (\ref{main}) written with the 
connection $\tilde \nabla.$ We have to prove that $Q^{\tilde \nabla}= 
Q^{\nabla}.$ 

We need some formul{\ae}. 

Since $\tilde \nabla$ is projectively equivalent to $\nabla$ there exists a 1-form $\omega$ on $M$ 
such that the Christoffel symbols of the connections $\nabla$ and $\tilde \nabla$ are related by 
\begin{equation}
\label{kob}
\tilde \Gamma_{ij}^k=\Gamma_{ij}^k+ \delta_{i}^k \omega_j + \delta_{j}^k \omega_i,
\end{equation}
(see \cite{kn}). It follows that, for any $\phi \in \cF_\l,$ one has $\nabla_k \phi= \tilde 
\nabla_k \phi + \l (1+n)\omega_k.$ In the same manner we can express the tensors 
$ \nabla_i\nabla_j \phi, \nabla_i\nabla_j\nabla_k \phi,$ with the tensors $\tilde \nabla_i
\tilde \nabla_j \phi, \tilde \nabla_i\tilde \nabla_j\tilde \nabla_k \phi,$ respectively. 

Using also formula (\ref{kob}) one has $\nabla_i P^{ijk}=\tilde \nabla_i P^{ijk}+ ((1+n)\delta-
(n+5))\omega_i P^{ijk}.$ In the same manner we can express the tensors $\nabla_j\nabla_i P^{ijk}, 
\nabla_k\nabla_j\nabla_i P^{ijk}, R_{ij}, \nabla_k R_{ij},$ with the tensors  
$\tilde \nabla_j\tilde \nabla_i P^{ijk}, \tilde \nabla_k\tilde \nabla_j\tilde 
\nabla_i P^{ijk},\tilde  R_{ij}, \tilde \nabla_k\tilde  R_{ij}.$ Replacing now these 
formul{\ae} into (\ref{main}), we see that $Q^{\nabla }=Q^{\tilde \nabla}$ if and only if the 
constants $\a, \beta_1, \beta_2, \eta_1,\eta_2,\eta_3$ are given as above. 

To prove part (ii), recall that the Lie algebra $\Sl_{n+1}(\bbR)$ can be identified with the Lie 
sub-algebra of $\Vect(\bbR^n)$ generated by the vector fields $\partial_i, x^i\partial_{j}, x^i x^j\partial _{j},$ 
where $(x^i)$ is the coordinates of the projective structure. The proof now is a simple computation 
(see \cite{lo}).\\
\cqfd

For the particular values of $\d:$
\begin{pro} If $\d=\frac{n+3}{n+1},\frac{n+4}{n+1},\frac{n+5}{n+1},$ there exists a quantization 
map given by (\ref{main}) with particular values of $\l,\mu,$ given in the following table
\end{pro}
\begin{center}
\label{tableau}
\begin{tabular}{ccccccccc}\hline 
$\d$ & $\l$ & $\mu$ \rule[-5mm]{0mm}{12mm}
 & $\alpha $  & $\b_1$  & $\b_2$ & $\eta_1$ & $\eta_2$ & $\eta_3$ \\\hline  
$\displaystyle \frac{n+5}{n+1}$ &  
$\displaystyle\frac{-2}{n+1}$ \rule[-5mm]{0mm}{12mm} & 
$\displaystyle\frac{n+3}{n+1}$ & $t$  & $ t $  &
$\displaystyle\frac{4}{1-n}$ &$\displaystyle \frac{1}{3}t $ &
$ \displaystyle \frac{4}{3}\frac{t }{(1-n)}$ & $\displaystyle 
\frac{2}{1-n}$  \\[2mm] 
$\displaystyle\frac{n+4}{n+1}$ & $\displaystyle \frac{-2}{n+1}$  &$\displaystyle\frac{n+2}{n+1}$
\rule[-5mm]{0mm}{12mm}   & $0$  &  $t $ & $\displaystyle \frac{4+t}{(1-n)}$  
 &$\displaystyle  \frac{2 }{3}t $ &$\displaystyle \frac{2}{3}\frac{ t}{(1-n)}$ & $\displaystyle 
\frac{6+2t}{(3-3n)}$ \\[2mm] 
$\displaystyle\frac{n+4}{n+1}$ &   $\displaystyle\frac{-1}{n+1}$ & 
\rule[-5mm]{0mm}{12mm}   $\displaystyle\frac{n+3}{n+1} $ & $3$ & $t$ & $\displaystyle\frac{1+t}{1-n}$& 
$\displaystyle\frac{1}{3}t$ &$\displaystyle\frac{9+t}{3-3n}$&$\displaystyle\frac{3+t}{3-3n}$\\[2mm] 
$\displaystyle\frac{n+3 }{n+1}$ &   $\displaystyle\frac{-2}{n+1}$ & 
\rule[-5mm]{0mm}{12mm}   $\displaystyle\frac{n+1}{n+1} $  &   $0$ &   $0$ &   $\displaystyle\frac{4}{1-n}$& 
$t$ & $\displaystyle \frac{4}{1-n}t$ &  $\displaystyle2\frac{1+t}{1-n}$
 \\[2mm] 
$\displaystyle\frac{n+3}{n+1}$ &   $\displaystyle\frac{-1}{n+1}$ & 
\rule[-5mm]{0mm}{12mm}   $\displaystyle\frac{n+2}{n+1} $  &   $\displaystyle \frac{3}{2}$ &   $0$ &  
 $\displaystyle\frac{1}{1-n}$& 
$t$ & $\displaystyle\frac{1}{2}\frac{(8t+3)}{(1-n)}$ &  $\displaystyle\frac{1+2t}{1-n}$
 \\[2mm] 
$\displaystyle\frac{n+3}{n+1}$ &   $0$ & 
\rule[-5mm]{0mm}{12mm}   $\displaystyle\frac{n+3}{n+1} $  &   $3$ &   $3$ &   $\displaystyle\frac{4}{1-n}$& 
$t$ & $\displaystyle\frac{4}{1-n}t$ &  $\displaystyle \frac{2}{1-n}t$
 \\[2mm] \hline
\end{tabular}
\end{center}
Here $t$ is a parameter.
\begin{rmk}{\rm (i) For the particular values of $\d,$ the quantization map (\ref{main}) is not unique 
(it is given by the parameter $t$). 

(ii) In the one dimensional case, the quantization map was given in \cite{cmz,ga}.

(iii) Another approach to the quantization map equivariant with respect to the action of the conformal group in a 
Riemannian manifold was given in \cite{do,l}.
}
\end{rmk}
\bigskip

{\it Acknowledgments}. I am grateful to Ch. Duval and V. Ovsienko for the statement of the problem. I am also 
grateful to H. Gargoubi and S. E. Loubon Djounga for fruitful discussions, and  Y. Maeda and Keio University for 
their hospitality. 

\vskip 1cm


\end{document}